\documentclass[reqno]{amsart}

\usepackage{amsmath,amsthm,amssymb,amsfonts}

\theoremstyle{plain}
\newtheorem{theo+}           {Theorem}      
\newtheorem{prop+}  [theo+]  {Proposition}
\newtheorem{coro+}  [theo+]  {Corollary}
\newtheorem{lemm+}  [theo+]  {Lemma}
\newtheorem{defi+}  [theo+]  {Definition}

\theoremstyle{definition}
\newtheorem{exam+}  [theo+]  {Example}
\newtheorem{rema+}  [theo+]  {Remark}

\newenvironment{theorem}{\begin{theo+}}{\end{theo+}}

\begin{document}

\baselineskip 18pt
\larger[2]
\title
[An elliptic determinant transformation]
{An elliptic determinant transformation} 
\author{Hjalmar Rosengren}
\address
{Department of Mathematical Sciences
\\ Chalmers University of Technology and G\"oteborg
 University\\SE-412~96 G\"oteborg, Sweden}
\email{hjalmar@math.chalmers.se}
\urladdr{http://www.math.chalmers.se/{\textasciitilde}hjalmar}

\thanks{Research  supported by the Swedish Science Research
Council (Vetenskapsr\aa det)}

\begin{abstract}
We prove a transformation formula relating two determinants involving
elliptic shifted factorials. Similar
determinants have been applied to multiple elliptic hypergeometric
series. 
\end{abstract}

\maketitle

\section{Introduction}  
Determinant evaluations play an important role in mathematics, perhaps
most notably in combinatorics, see Krattenthaler's surveys \cite{K} 
and \cite{K2}. Many useful determinant evaluations
are \emph{rational} identities, which rises the question of finding
 generalizations to the \emph{elliptic} 
level. In recent work with Schlosser \cite{RS}, we gave an
approach to elliptic determinant evaluations that encompasses most results in the literature, 
from the classical Frobenius determinant  to the
Macdonald identities  for  non-exceptional affine root systems. 

As an example of an elliptic determinant evaluation, we mention
Warnaar's determinant \cite[Corollary 5.4]{W}, which we write as
\begin{multline}\label{wd}\det_{1\leq j,k\leq
  n}\left(\frac{(bx_j,c/x_j)_{k-1}}{(a/bx_j,ax_j/c)_{k-1}}\right)\\
  =c^{\binom n2}q^{\binom n3}\prod_{1\leq i<j\leq n}
x_i^{-1}\theta(x_i/x_j)\theta(bx_ix_j/c)\prod_{j=1}^n\frac{(a/bc,aq^{j-2})_{j-1}}{(a/bx_j,ax_j/c)_{n-1}}.
 \end{multline}
Here, we use the  notation
$$\theta(x)=\prod_{j=0}^\infty(1-p^jx)(1-p^{j+1}/x), $$
$$(a)_k=\theta(a)\theta(aq)\dotsm\theta(aq^{k-1}),$$
$$(a_1,\dots,a_n)_k=(a_1)_k\dotsm(a_n)_k, $$
where  $p$ and $q$ are fixed parameters with $|p|<1$.
In the trigonometric case, $p=0$,  we recover the usual $q$-shifted
factorials, which we  denote
$$(a)_k^{\text{trig}}=(1-a)(1-aq)\dotsm(1-aq^{k-1}). $$
The case $p=0$ of \eqref{wd} is a special case of a determinant
evaluation due to Krattenthaler \cite[Lemma~34]{Kr}.

Warnaar used \eqref{wd} to derive a  multivariable extension of
the elliptic Jackson summation.  In the terminology of
\cite{DS}, this is a ``Schlosser-type'' sum, which is obtained by taking
the determinant of one-dimensional summations. In spite of its
conceptual simplicity, Warnaar's sum is a key result for
multiple elliptic hypergeometric series, since it can be used both to derive
 an ``Aomoto--It\^o--Macdonald-type'' sum (conjectured in \cite{W} and proved in \cite{R1}) and a ``Gustafson--Milne-type''
sum (conjectured in \cite{DS} and  proved
in \cite{R2}). Alternative proofs of these summations were found
 by Rains \cite{Ra1,Ra2}. 
 See  \cite{Sp} for an application of
\eqref{wd} to elliptic hypergeometric integrals.

In \cite[Eq.~(7.27)]{RS0}, as a by-product of deriving 
transformation formulas for Schlosser-type  series, we discovered
 the identity
\begin{equation}\label{tdt}
\det_{1\leq j,k\leq n}\left(\frac{(z_j)_{k-1}^{\text{trig}}}{(a_jz_j)_{k-1}^{\text{trig}}}\right)
=(-1)^{\binom n2}q^{\binom n3}\det_{1\leq j,k\leq
  n}\left(z_j^{k-1}\frac{(a_j)_{k-1}^{\text{trig}}}{(a_jz_j)_{k-1}^{\text{trig}}}\right). 
 \end{equation}
 Krattenthaler found a more natural proof of \eqref{tdt}, based on the
$q$-Chu--Vandermonde summation, which was included in \cite{RS0}.

The purpose of the present note is to obtain an
 elliptic extension of \eqref{tdt}. 
In fact, as we will explain at the end, such an identity can be derived from the results of \cite{RS0}. However,
we prefer to give a self-contained proof, 
which is a straight-forward
 generalization of Krattenthaler's proof of \eqref{tdt}, based on the
 elliptic Jackson summation
\begin{equation}\label{js}
\sum_{l=0}^n\frac{\theta(aq^{2l})}{\theta(a)}\frac{(a,b,c,d,e,q^{-n})_l}
{(q,aq/b,aq/c,aq/d,aq/e,aq^{n+1})_l}=\frac{(aq,aq/bc,aq/bd,aq/cd)_n}
{(aq/b,aq/c,aq/d,aq/bcd)_n},
\end{equation}
where $a^2q^{n+1}=bcde$. The identity \eqref{js} was obtained by Date et al.\ \cite{D} for
special parameter values and by Frenkel and Turaev \cite{FT} in
general;  see \cite{R3} for an elementary proof.

\begin{theorem}\label{th} Let $a$ and $b_j,c_j,d_j$, $j=1,\dots,n$, be
 parameters such that the product $b_jc_jd_j$ is independent of $j$. Then 
the  determinant transformation
\begin{multline}\label{dt}
\det_{1\leq j,k\leq n}\left(\frac{(b_j,c_j,d_j)_{k-1}}{(a/b_j,a/c_j,a/d_j)_{k-1}} \right)\\
=\left(\frac{a}{e}\right)^{\binom n2}\prod_{j=2}^n
\frac{(aq^{j-2})_{j-1}}{(eq^{j-2})_{j-1}}
\det_{1\leq j,k\leq n}\left(\frac{(a/b_jc_j,a/b_jd_j,a/c_jd_j)_{k-1}} 
{(a/b_j,a/c_j,a/d_j)_{k-1}}\right)
\end{multline}
holds, where 
\begin{equation}\label{bal}e=a^2/b_jc_jd_j.\end{equation}
\end{theorem}

Note that if $d_j$ is independent of $j$, we may write $b_j=bx_j$,
$c_j=c/x_j$ and evaluate both determinants  using \eqref{wd}.
Thus, in this special case, 
Theorem \ref{th}
follows from Warnaar's determinant evaluation.

As far as we know, Theorem \ref{th} is new even in the case $p=0$. In
that case,   letting $a\rightarrow 0$,
$c_j\rightarrow 0$ and $d_j\rightarrow\infty$,
keeping $b_j$ and $a/c_j$ fixed, so that $e\rightarrow 0$,
 one recovers \eqref{tdt} after a change
of variables.

In view of its close relation to Warnaar's determinant, one may hope that
 Theorem \ref{th} will 
also find applications to multiple elliptic hypergeometric
series. However, so far we have not found any interesting results in
that direction. Perhaps Theorem \ref{th} serves to indicate that when
encountering a determinant that cannot be evaluated in closed form,
one should keep in mind that
it may still satisfy some,  potentially useful, transformation.

\section{Proof of Theorem \ref{th}}

Let 
$X=(X_{jk})_{j,k=1}^n$ and $Y=(Y_{jk})_{j,k=1}^n$ be the matrices
\begin{equation}\label{x}X_{jk}=\frac{(b_j,c_j,d_j)_{k-1}}{(a/b_j,a/c_j,a/d_j)_{k-1}},\end{equation}
$$Y_{jk}=\frac{\theta(aq^{2j-3})}{\theta(aq^{-1})}\frac{(aq^{-1},q^{1-k},eq^{k-2})_{j-1}}{(q,aq^{k-1},aq^{2-k}/e)_{j-1}}\,q^{j-1}. $$
Note that $Y$ is triangular, with determinant
\begin{equation}\label{ydet}\begin{split}\det(Y)&=\prod_{j=1}^n
  Y_{jj}=\prod_{j=2}^n\frac{\theta(aq^{2j-3})}{\theta(aq^{-1})}\frac{(aq^{-1},q^{1-j},eq^{j-2})_{j-1}}{(q,aq^{j-1},aq^{2-j}/e)_{j-1}}\,q^{j-1}\\
&=\left(\frac{e}{a}\right)^{\binom n2}\prod_{j=2}^n
\frac{(a,eq^{j-2})_{j-1}}{(aq^{j-2},e/a)_{j-1}},
 \end{split}\end{equation}
where we used the elementary identity
\begin{equation}\label{ei}\frac{(x)_{j-1}}{(y)_{j-1}}=\left(\frac xy\right)^{j-1}\frac{(q^{2-j}/x)_{j-1}}{(q^{2-j}/y)_{j-1}}\end{equation}
in the last step. Moreover,
$$(XY)_{jk}=\sum_{l=0}^{n-1}
X_{j,l+1}Y_{l+1,k}=\sum_{l=0}^{k-1}
\frac{\theta(aq^{2l-1})}{\theta(aq^{-1})}
\frac{(aq^{-1},q^{1-k},b_j,c_j,d_j,eq^{k-2})_{l}}{(q,aq^{k-1},a/b_j,a/c_j,a/d_j,aq^{2-k}/e)_{l}}\,q^{l}.
 $$
Assuming \eqref{bal}, the elliptic Jackson summation \eqref{js} gives
\begin{equation}\label{xy}(XY)_{jk}=
\frac{(a,a/b_jc_j,a/b_jd_j,a/c_jd_j)_{k-1}} 
{(e/a,a/b_j,a/c_j,a/d_j)_{k-1}}.\end{equation} 
Writing out the equation $\det(X)=\det(XY)/\det(Y)$ using \eqref{x},
\eqref{ydet} and \eqref{xy} we arrive at \eqref{dt}.

\section{An $S_3$ symmetry}

Besides the non-trivial symmetry of Theorem \ref{th}, determinants of
the form
$$\det_{1\leq j,k\leq
    n}\left(\frac{(b_j,c_j,d_j)_{k-1}}{(a/b_j,a/c_j,a/d_j)_{k-1}}
  \right),\qquad  b_jc_jd_j \text{ independent of $j$} ,
$$
 also have a trivial symmetry. Namely, reversing the
order of the columns, we have
\begin{multline*}\det_{1\leq j,k\leq
    n}\left(\frac{(b_j,c_j,d_j)_{k-1}}{(a/b_j,a/c_j,a/d_j)_{k-1}}
  \right)=(-1)^{\binom n2} \det_{1\leq j,k\leq
    n}\left(\frac{(b_j,c_j,d_j)_{n-k}}{(a/b_j,a/c_j,a/d_j)_{n-k}}\right)\\
=(-1)^{\binom
  n2}\prod_{j=1}^n\frac{(b_j,c_j,d_j)_{n-1}}{(a/b_j,a/c_j,a/d_j)_{n-1}}\\
\times \det_{1\leq j,k\leq
    n}\left(\frac{(q^{n-k}a/b_j,q^{n-k}a/c_j,q^{n-k}a/d_j)_{k-1}}{(q^{n-k}b_j,q^{n-k}c_j,q^{n-k}d_j)_{k-1}}\right).
\end{multline*}
Using \eqref{ei}
and introducing the parameter $e$ as in \eqref{bal} gives
\begin{multline}\label{ts}\det_{1\leq j,k\leq
    n}\left(\frac{(b_j,c_j,d_j)_{k-1}}{(a/b_j,a/c_j,a/d_j)_{k-1}}
  \right)
=\left(-\frac{e^2}{a}\right)^{\binom
  n2}\prod_{j=1}^n\frac{(b_j,c_j,d_j)_{n-1}}{(a/b_j,a/c_j,a/d_j)_{n-1}}\\
\times\det_{1\leq j,k\leq
    n}\left(\frac{(q^{2-n}b_j/a,q^{2-n}c_j/a,q^{2-n}d_j/a)_{k-1}}{(q^{2-n}/b_j,q^{2-n}/c_j,q^{2-n}/d_j)_{k-1}}
  \right).
 \end{multline}

Denoting by $\sigma$ and $\tau$ the transformation from the
left-hand to the right-hand side of \eqref{dt} and \eqref{ts},
respectively,  one may check that 
$\sigma^2=\tau^2=(\sigma\tau)^3=\operatorname{id}$, that is,
$\sigma$ and $\tau$ generate an  $S_3$ symmetry. 
 Thus, there are  three additional expressions,
 corresponding to $\sigma\tau$, $\tau\sigma$ and
 $\sigma\tau\sigma=\tau\sigma\tau$. These may be written,
 respectively, as
\begin{multline}\det_{1\leq j,k\leq
    n}\left(\frac{(b_j,c_j,d_j)_{k-1}}{(a/b_j,a/c_j,a/d_j)_{k-1}}
  \right)\\
\begin{split}&=q^{-6\binom n3}\left(\frac {e}{a^2}\right)^{\binom n2}
\prod_{j=1}^n\frac{(b_j,c_j,d_j)_{n-1}}{(a/b_j,a/c_j,a/d_j)_{n-1}}
\prod_{j=2}^n\frac{(aq^{j-2})_{j-1}}{(q^{j-n}e/a)_{j-1}}\\
&\quad\times\det_{1\leq j,k\leq
    n}\left(\frac{(a/b_jc_j,a/b_jd_j,a/c_jd_j)_{k-1}}{(q^{2-n}/b_j,q^{2-n}/c_j,q^{2-n}/d_j)_{k-1}}\right)\\
&=\left(-\frac {a^3}{e^2}\right)^{\binom n2}
\prod_{j=1}^n\frac{(a/b_jc_j,a/b_jd_j,a/c_jd_j)_{n-1}}{(a/b_j,a/c_j,a/d_j)_{n-1}}
\prod_{j=2}^n\frac{(aq^{j-2})_{j-1}}{(eq^{j-2})_{j-1}}\\
&\quad\times\det_{1\leq j,k\leq
    n}\left(\frac{(q^{2-n}b_j/a,q^{2-n}c_j/a,q^{2-n}d_j/a)_{k-1}}{(q^{2-n}b_jc_j/a,q^{2-n}b_jd_j/a,q^{2-n}c_jd_j/a)_{k-1}}\right)\\
&=q^{-6\binom n3}\left(\frac {a^2}{e^3}\right)^{\binom n2}
\prod_{j=1}^n\frac{(a/b_jc_j,a/b_jd_j,a/c_jd_j)_{n-1}}{(a/b_j,a/c_j,a/d_j)_{n-1}}
\prod_{j=2}^n\frac{(aq^{j-2})_{j-1}}{(q^{j-n}a/e)_{j-1}}\\
&\quad\times\det_{1\leq j,k\leq
    n}\left(\frac{(b_j,c_j,d_j)_{k-1}}{(q^{2-n}b_jc_j/a,q^{2-n}b_jd_j/a,q^{2-n}c_jd_j/a)_{k-1}}\right).\\
\end{split}\label{et}
\end{multline}
 To verify this, the elementary
identities
$$\prod_{j=2}^n(aq^{j-2})_{j-1}=\prod_{j=2}^n(aq^{2n-2j})_{j-1}
=(-a)^{\binom n2}q^{3\binom n3}
\prod_{j=2}^n(q^{2-2n+j}/a)_{j-1}$$
are useful. Although each of 
the three transformations in \eqref{et} are  equivalent
to Theorem~\ref{th},  
it may be worthwhile to state them explicitly.

Finally, we explain how to obtain Theorem \ref{th} from the results of
\cite{RS0}. 
In that paper (the elliptic extension of Corollary 7.3 contained in Theorem 8.1), we obtained the identity 
\begin{multline}\label{cnt}
\sum_{k_1,\dots,k_n=0}^{m_1,\dots,m_n}\,
\prod_{1\le i<j\le n}q^{k_i}\theta(q^{k_j-k_i})
  \theta(aq^{k_i+k_j})\\
\times \prod_{j=1}^n\frac{\theta(aq^{2k_j})}{\theta(a)}\frac{(a,b,c_j,d_j,e_j,q^{-m_j})_{k_j}}
{(q,aq/b,aq/c_j,aq/d_j,aq/e_j,aq^{1+m_j})_{k_j}}\,
q^{k_j}\\
=b^{-\binom n2}q^{-2\binom n 3}
\prod_{j=1}^n\frac{(aq^{2-n}/b)_{n-1}(b)_{j-1}}{(aq^{2+n-2j}/b)_{n-1}}
\frac{(aq,aq/c_jd_j,aq/c_je_j,aq/d_je_j)_{m_j}}
{(aq/c_j,aq/d_j,aq/e_j,aq/c_jd_je_j)_{m_j}}\\
\times\det_{1\le j,k\le n}\left(
\frac{(c_j,d_j,e_j,q^{-m_j})_{k-1}}
{(aq^{2-n}/bc_j,aq^{2-n}/bd_j,aq^{2-n}/be_j,aq^{2-n+m_j}/b)_{k-1}}\right),
\end{multline}
where $bc_jd_je_j=a^2q^{2-n+m_j}$ for $j=1,\dots,n$. Note that 
the case $n=1$ is \eqref{js}. Consider the special case 
 when $m_i=n-1$ for each $i$. Then, the terms in the sum
vanish unless  $(k_1,\dots,k_n)$
is a permutation of $(0,1,\dots,n-1)$. Moreover, since
$$\prod_{1\le i<j\le n}q^{k_i}\theta(q^{k_j-k_i})=
\operatorname{sgn}(k)\prod_{0\le i<j\le n-1}q^i\theta(q^{j-i}),$$
the sum can be written as a determinant. After some elementary computation
and a change of variables, one is reduced to  the 
equality between the first and last member of \eqref{et}. Thus,
Theorem \ref{th} can also be obtained as a special case of \eqref{cnt}.

\end{document}